\documentclass[12pt]{article}
\usepackage{theorem}
\usepackage{amssymb}
\usepackage{amsmath}
\usepackage{epsfig}
\usepackage{graphicx}

\newtheorem{theo}{Theorem}[section]
\newtheorem{prop}[theo]{Proposition}
\newtheorem{coro}[theo]{Corollary}
\newtheorem{lem}[theo]{Lemma}

{\theorembodyfont{\rm}

\newtheorem{dfn}[theo]{Definition}

\newtheorem{rem}[theo]{Remark}

}

\def\ra{\rightarrow}

\def\C{{\mathbb C}}
\def\S{{\mathbb S}}
\def\D{{\Delta}}
\def\E{{\cal E}}
\def\F{{\mathbb F}}
\def\G{{\Gamma}}
\def\Ga{{\Gamma}}
\def\N{{\mathbb N}}

\def\P{{\mathbb P}}
\def\Q{{\mathbb Q}}
\def\R{{\mathbb R}}

\def\Z{{\mathbb Z}}

\def\Hom{{\rm Hom}}

\def\SL{{\rm SL}}
\def\PSL{{\rm PSL}}

\def\H{{\cal H}}
\def\no{\noindent}
\def\deg{{\rm deg}}
\def\df{{\rm DF}}
\def\cdf{{\rm CDF}}
\def\Ker{{\rm Ker}}
\def\ET{{\rm ET}}

\author{Fedor Bogomolov and Yuri Tschinkel}
\title{{\Large Monodromy of elliptic surfaces}}
 
\begin{document}
 
\maketitle
 
\section{Introduction}
\label{sect:intro}

Let ${\E}\ra B$ be a non-isotrivial 
Jacobian elliptic fibration 
and $\tilde{\Gamma}$ its global monodromy group.
It is a subgroup of finite index in $ \SL(2,\Z)$. 
We will assume that ${\E}$ is Jacobian. Denote by $\Ga$ the
image of $\tilde{\Ga}$ in $\PSL(2,\Z)$ and by
$\overline{\cal H}$ the upper half-plane completed by $\infty$ and by
rational points in $\R\subset \C$. 
The $j$-map $B\ra \P^1$ decomposes as 
$j_{\Ga}\circ j_{\E}$,
where 
$$
j_{\cal E} \,:\, B\ra M_{\Ga}=\overline{\cal H}/{\Ga}
$$ and
$j_{\Gamma}\,:\, M_{\Gamma} \ra \P^1= 
\overline{\cal H}/\PSL(2,\Z)$. 
In an algebraic family of elliptic fibrations the degree 
of $j$ is bounded by the degree of the generic element.
It follows that there is only a finite number 
of monodromy groups for each family.  
 
The number of subgroups of bounded index in $\SL(2,\Z)$ 
grows superexponentially \cite{lubotzky}. However, 
the number of $M_{\Gamma}$-representations of the sphere $\S^2$
grows exponentially. Thus, monodromy groups of elliptic
fibrations over $\P^1$ constitute a small, but still
very significant fraction of all subgroups 
of finite index in $\SL(2,\Z)$.

Our goal is to introduce some structure
on the set of monodromy groups of elliptic
fibrations which would help to answer some natural 
questions.  
For example, we show how to describe the set
of groups corresponding to rational or K3 elliptic
surfaces, explain  how to compute the dimensions of 
the spaces of moduli of surfaces in this class with 
given monodromy group etc. 
Our method is based on 
a detailed study of triangulations of Riemann 
surfaces.

To determine $\tilde{\Gamma}$ we first describe all possible 
groups ${\Gamma}$. 
In order to classify possible ${\Gamma}$
we consider the corresponding oriented Riemann surface 
$M_{\Gamma}$. The map $j_{\Gamma}\,:\, M_{\Gamma}\ra \P^1$
provides a special triangulation of $M_{\Gamma}$ 
(induced from the standard triangulation of $\P^1$ into
two triangles with vertices in $0,1,\infty$) (and vice versa).  
The preimages of $0,1,\infty$ on 
$M_{\Gamma}$ will be called $A,B,I$, respectively,
and the triangulation will be called an 
$ABI$-triangulation.    
The barycentric subdivision of any triangulation of 
an oriented Riemann surface is an $ABI$-triangulation. 
This remark goes back at least to Alexander \cite{alexander}
(who proves an analogous statement in {\em any} dimension). 
Constructions of this type were 
rediscovered by many authors in connection with Belyi's theorem
and Grothendieck's ``Dessins d'enfants'' program
(\cite{belyi}, \cite{schneps-lochak} and the references therein). 
An $ABI$-triangulation of a Riemann surface $R$ induces
a graph on $R$,  
which is obtained by removing all  
$AI$- and $BI$-edges  from the graph given by the 1-skeleton 
of the $ABI$-triangulation
(see \cite{shabat-voevodski}, for 
example). In our case, we have
additional constrains on the valence of the $A,B$ vertices in this graph.
Namely, the $A$-vertices have valence 1 or 3,
and the  $B$-vertices have valence 1 or 2. More constrains arise
from considerations of local monodromies.

The plan of the paper is as follows.
In section \ref{sect:generalities} we recall basic facts about
the local and global monodromy groups of elliptic fibrations
due to Kodaira. In section \ref{sect:modular-curves} we study
$j$-modular curves $M_{\Ga}$ and their relationship with
$ABI$-triangulations. 
In section \ref{sect:modular-surfaces}
we give a modular construction of elliptic surfaces over $M_{\Ga}$
with prescribed monodromy groups. 
General elliptic fibrations over $B$ are obtained
as simple modifications of pullbacks of these elliptic fibrations
from $M_{\Ga}$ - this is the content of section \ref{sect:lifts}. 
Our construction allows a relatively transparent
description of a rather complicated set of global  
monodromy groups of elliptic surfaces. This transforms the general
results of Kodaira theory to a concrete computational tool.

\

\no
{\bf Conventions.} 
We write $\Z_n=\Z/n\Z$ and $\F_n$ for the free group on $n$
generators.  Throughout the paper we work over $\C$.

\

\no
{\bf Acknowledgments.} The first author was partially 
supported by the NSF. The second author was partially
supported by the NSA. We are grateful to J. de Jong and 
R. Vakil for helpful discussions.

\section{Generalities}\label{sect:generalities}

In this section we give a brief summary of Kodaira's theory
of elliptic fibrations. 
We refer to the papers by Kodaira \cite{kodaira} and to 
\cite{barth-peters-vdv} and \cite{isk-sha} 
for proofs and details. 

\subsection{The setup}
Let 
$f \,:\, \E\ra B$  be a smooth 
relatively minimal non-isotrivial Jacobian 
elliptic fibration over a smooth curve $B$ of genus $g(B)$. 
This means that 
\begin{enumerate}
\item  $\E$ is a smooth compact surface and $f$ is holomorphic,
\item the generic fiber of $f$ is a 
smooth curve of genus 1 (elliptic fibration),
\item the fibers of $\E$ do not contain
smooth rational curves of self-intersection $-1$ (relative minimality),
\item  we have a global zero section 
$e\, :\,B\ra \E$ (Jacobian elliptic fibration),
\item the $j$-function which to  each smooth  fiber $\E_b\subset {\cal E}$
assigns its $j$-invariant is a non-constant 
rational function on $B$ (non-isotrivial). 
\end{enumerate}

\subsection{Topology}

Denote by $B^s=\{b_1,...,b_k\}\subset B$ the set of points
corresponding to singular fibers of $\E$, it is always non-empty. 
Let $B^0=B\backslash B^s$ be the open 
subset of $B$ where all fibers are smooth
and $f^0\,:\, \E^0\ra B^0$ the restriction of $f$. 
Topologically, $f^0$ is a smooth oriented 
fibration with fibers $\S^1\times \S^1$, 
which is equipped with a section. 
The equivalence class of 
$\E^0$ under global diffeomorphisms inducing smooth
isomorphisms on each fiber is 
determined by the topology of $B^0$ and by the
homomorphism (representation) 
of the fundamental group $\pi_1(B^0)$ into
the group of  homotopy classes of  
orientation-preserving automorphisms of the torus
$\S^1\times \S^1$. The latter is isomorphic  to $\SL(2,\Z)$ and 
this isomorphism is uniquely determined 
by the choice of generators of 
$\pi_1(\S^1\times \S^1)=\Z\oplus \Z$ (with fixed orientation).  
A different choice of generators leads to a conjugation (in 
$\SL(2,\Z)$) of the isomorphism.
Thus we have a homomorphism 
$\rho_{\E}^c\,:\, \pi_1(B^0)\ra \SL(2,\Z)$. 
This homomorphism - it is defined 
modulo conjugation in $\SL(2,\Z)$ - is called
by Kodaira the {\em homological invariant} 
of the elliptic fibration $\E$.

Now we consider the local situation: according to Kodaira, 
the restriction of $f$ to a 
small punctured analytic neighborhood $\D_b^*$
of a point $b\in B$ (disc $\D_b$ minus the  point $b$) for 
every point $b\in B^s$ is also topologically non-trivial. 
Thus we have a homomorphism 
$\rho_b^c\,:\, \Z\ra \SL(2,\Z)$  (where $\Z$ is the 
fundamental group of the punctured disc 
with the standard generator $t_b$). 
Again, this homomorphism is defined modulo
conjugation.

We can eliminate the ambiguity in the definitions above
by the following procedure: choose a point $b_0\in  B^0$ 
and a set of non-intersecting
paths   connecting $b_0$ to the singular points $b_s\in B^s$. 
This set admits a natural cyclic order defined by the 
relative position of  these paths
in a small neighborhood of $b_0$.

A small  neighborhood of this set is a 
disc inside $B$ (with orientation). Now   we can 
choose small oriented loops 
around each singular point $b_s$.

\vskip 0,5cm
\hskip 3,5cm\includegraphics{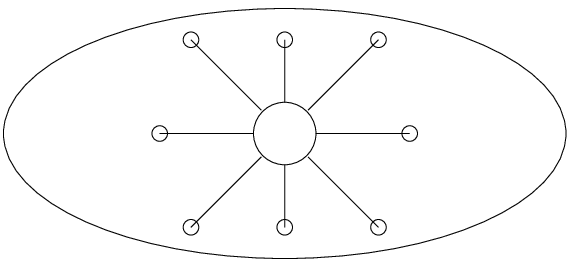}
\vskip 0,5cm

If we fix generators  of $\pi_1(\S^1\times \S^1)$ 
for the fiber over $b_0$ then we
obtain a system of {\em elements} $T_b\in \SL(2,\Z)$ 
in the conjugacy class of $t_b$ as well as a representation 
$\rho_{\E} \,:\, \pi_1(B^0)\ra \SL(2,\Z)$. 
We call the elements $T_b$ local monodromies, the representation 
$\rho_{\E}$ the global monodromy representation and the 
group 
$\tilde{\Ga}=\rho_{\E}(\pi_1(B^0))\subset \SL(2,\Z)$ 
the global monodromy group. 
The global monodromy representation depends 
only on the basis of $\pi_1(\S^1\times \S^1)$
at $b_0$. The local monodromy elements depend 
on the choice of the system of paths.

There is an important relation between local and global monodromy.

\begin{lem} Let $\E\ra B$ be an elliptic fibration as above. 
Suppose that $B=\P^1$. Then the product  
$$
P(\E):=\prod_{b\in B^s}T_b \in \SL(2,\Z)
$$ 
(taken in cyclic order) is equal to the identity. 
Similarly, if the genus  $g(B)\ge 1$ 
then  $P(\E)$ is a product of $g(B)$ commutators.  
\end{lem}

{\em Proof.}
The product $P(\E)$ gives the 
monodromy transformation along the boundary
of the disc $\D$. Our fibration is smooth on the 
complement $B\setminus \D$. Therefore, it
is a topologically trivial fibration over a 
disc in the case of $B=\P^1$ or a smooth 
$\S^1\times \S^1$ fibration over the Riemannian surface 
$B$ of genus $g(B)$ minus a
disc. Now the relations follow from similar 
relations in $\pi_1(B\setminus \D)$.

\subsection{The $j$-function}
\label{sect:j-map}


The elliptic fibration $\E\ra B$ defines a rational 
function on $B$ - the $j$-function. 
There is a relationship between the $j$-function 
and local (resp. global) monodromies. 

First we look at the local situation: the
restriction of $j$ to the disc 
$\D_b$ is analytically equivalent to 
$j(b)+z^k$  if $j(b)$ is finite
or $z^{-k}$ if $j(b)$ is infinite ($k\in \N$). 
Here $z$ is a local parameter.
There are certain compatibility conditions between $k$ 
and the local monodromy $\rho_b$. 
Kodaira gives a list of all pairs $(\rho_b,k)$ 
which occur. Moreover,  

\begin{theo}
The pair $(\rho_b, k)$ defines a unique (in the
analytic category) semistable Jacobian 
fibration over $\D_b$. 
Any two Jacobian elliptic fibrations 
over an analytic disc $\D_b$ 
with the same $(k, \rho_b)$ 
are fiberwise birationally isomorphic. 
\end{theo}

Globally, the $j$-function determines the {\em image}
of the global monodromy $\rho_{\E}^c$
in $\PSL(2,\Z)$. 
There are exactly $2^{g(B)+k-1} $ different liftings of 
the standard generators of $\pi_1(B^0) $  to $\SL(2,\Z)$,
which correspond to homomorphisms of $\pi_1(B^0)$ into $\SL(2,\Z)$.
The local liftings differ by the central 
element $c\in \Z_2\subset \SL(2,\Z)$. 
Each of these liftings determines a unique homological invariant,
admissible for $j$. All of $j$-admissible homological invariants 
are obtained in this way. This explains the part (a) of the theorem 
11.1 p. 160 in \cite{barth-peters-vdv}.

\begin{theo}
\label{theo:global}
Let $B$ be a connected compact curve and $b_1,...,b_k$ a finite
set of points on $B$. Let $j$ be a non-constant rational function on 
$B$ such that $j\neq 0,1,\infty$ on 
$B^0= B\setminus \{ b_1,...,b_k\}$.
For a fixed homological invariant $\rho_{\E}^c$, 
which is admissible for $j$, there exists a unique 
Jacobian elliptic fibration $\E$ with 
this $j$ and $\rho_{\E}^c$. 
\end{theo}

Suppose however, that we are interested in classifying global
monodromies in some restricted class of surfaces, for example 
rational elliptic or elliptic K3 surfaces. Then only a finite
number of possible global monodromy groups $\tilde{\Ga}$ 
and only few homological invariants can occur if we fix the image
of $\tilde{\Ga}$  in $\PSL(2,\Z)$.
It is clear that elliptic surfaces with the same $j$-invariant 
but different homological invariants are scattered through different
topological classes. 
Our point of departure was that in this situation the 
Theorem \ref{theo:global} does not provide any simple and 
sufficient control over the topology of the resulting surfaces. 
In the following sections we give some technical improvements
of Kodaira's theory which lead to an effective algorithm.

\section{$j$-modular curves}\label{sect:modular-curves}

To an  elliptic fibration $f\,:\,\E\ra B$ we can associate a 
curve defined over  $\overline{\Q}$ equipped with a special
triangulation. This triangulation will be our main tool in the
description of global monodromy $\tilde{\Ga}$ of  $\E$. 

Let $\Gamma$ be the image of $\tilde{\Gamma}$ in 
$\PSL(2,\Z) = \SL(2,\Z)/\Z_2$.
We have the $j$-map 
$ j : B \to \overline{\H}/\PSL(2,\Z) = \P^1$. 
This map decomposes as a product
$ j = j_{\Gamma}\circ j_{\E} $ where 
$j_{\E} : B\to \overline{\H}/\Gamma$ is a natural lifting
of $j$ onto the modular curve 
$M_{\Gamma} = \overline{\H}/\Gamma $    
corresponding to $\Gamma$ and
\begin{equation}\label{eqn:decomp}
j_{\Gamma} : \overline{\H}/\Gamma 
\ra \overline{\H}/\PSL(2,\Z) = \P^1.
\end{equation}
The above decomposition shows 
that $\deg(j) = \deg(j_{\E})\cdot \deg(j_{\Gamma})$. 
In particular, for 
any non-isotrivial elliptic surface the group $\Gamma$ 
is a subgroup of finite index in $\PSL(2,\Z)$.

\begin{dfn}
We call the pair $(M_{\Gamma}, j_{\Ga})$  the
$j$-modular curve corresponding to the
monodromy group $\Gamma$.
\end{dfn}

\begin{rem}
Usual modular curves are $j-modular$. 
A $j-modular$
curve is simply any curve defined over a number field together
with a special rational function on it
(this follows from the theorem of 
Belyi \cite{belyi}, see \ref{lem:jmodular}). 
There is a countable number of such functions for each curve. 
\end{rem}

Let us give a combinatorial 
description of $j$-modular curves.
They correspond to special 
triangulations of Riemann surfaces.

\begin{dfn}
Let $R$ be an oriented Riemann surface. 
A triangulation  $\tau(R)=(\tau_0,\tau_1,\tau_2)$ of $R$ is 
a decomposition of $R$ into a finite union 
of open 2-cells $\tau_2$ and a connected 
graph $\tau_1$ with vertices $\tau_0$ such that the complement
$\tau_1\setminus \tau_0$ is a disjoint union of open
segments and the closure of any open 2-cell is isomorphic to 
the image of a triangle under a simplicial map. 
\end{dfn} 

The number of edges originating in a vertex $x$ is
called the valence at $x$ and is denoted by $v(x)$.

\begin{dfn}
An $ABI$-triangulation of $R$ is a triangulation 
together with a coloring of vertices  
in three colors $A,B$ and $I$ such that
\begin{enumerate} 
\item The colors of any two adjacent vertices are different. 
\item There are $2$ or $6$ edges at vertices
of color $A$ and $2$ or $4$ edges at vertices of color $B$.
\end{enumerate}
\end{dfn}

We will refer to vertices of color $A$
(resp. $B$) with valence $j$ as $A_j$ (resp. $B_j$) vertices.  
If we delete the $I$-vertices from $\tau_0$ and 
all edges $AI$ and $BI$ from 
$\tau_1$ then the remaining  connected graph on $R$ is called the 
$AB$-graph associated to the $ABI$-triangulation. The valences
of $A$-vertices in an $AB$-graph are $1$ or $3$, the valences of
$B$-vertices are $1$ or $2$ and vertices of the same color are not 
connected by an edge.  The $I$-vertices from the $ABI$-triangulation 
are represented by connected components of $R$ minus the $AB$-graph.
An $ABI$-triangulation on $R$ can be 
reconstructed from an $AB$-graph on 
$R$ by placing one $I$-vertex into 
each connected component of $R$ minus
the $AB$-graph and by connecting  (cyclically) the 
$I$-vertex with vertices on the
boundary of the corresponding connected component.  

The following well known theorem forms the basis
for our analysis of monodromy groups. 

\begin{theo}\label{theo:modular}
Let $R$ be an oriented compact Riemann surface with
an $ABI$-triangulation.
Then there exists a unique structure of a
$j$-modular curve on $R$. Conversely, 
every structure of a $j$-modular curve on $R$
corresponds to an $ABI$-triangulation.
\end{theo}

{\em Proof.} 
Let us first show how $j_{\Ga}$ defines a triangulation of
$M_{\Gamma}$.
The map $j: \overline{\H}\to \overline{\H}/\PSL(2,\Z) = \P^1$ 
is ramified over
three points $ 0 = A, 1 = B, \infty = I$.
The ramification index at $0$ is equal to $2$, the ramification
index at $1$ is $3$ and the 
ramification index at $\infty$ is infinite. 
Similar result is true for
$$
j_{\Ga} : \overline{\H}/\Gamma =  M_{\Gamma} 
\to \overline{\H}/\PSL(2,\Z) = \P^1.
$$
Consider the standard triangulation $\tau_{st}(\S^2)$
of the sphere $\S^2 = \P^1$ into a union of
two triangles with vertices $0,1$ and $\infty$.
The preimage of this triangulation provides a  
triangulation of $M_{\Gamma}$. 
If we color the preimages of the
corresponding vertices 
in $A, B$ and $I$ 
then we obtain an $ABI$-triangulation as wanted.

Conversely, starting with an $ABI$-triangulation $\tau$
we construct an algebraic curve
$R$ together with a 
map $R\ra \S^2$ ramified in $0,1,\infty$ as follows.
We have  a map from the set of vertices to $ (A,B,I)$ (the color). 
Further, every edge will be mapped into the edges of 
the standard triangulation of $\S^2$, respecting the colors of the
ends. This map is completed by the map of triangles,
which maps the triangles $ABI$ (with orientation inherited
from $R$) to one of the triangles of $\tau_{st}(\S^2)$ and the
triangles with the opposite $R$-orientation 
to the other. 

Thus we have constructed a 
simplicial map  which is locally
an isomorphism except in the neighborhood of vertices.
Since triangles in  $R$ sharing an edge are mapped into different
triangles of $\S^2$ the above map is locally 
an isomorphism outside of vertices and is 
equivalent to a map $z^n$ in the neighborhood
of each vertex in $R$.
Thus it corresponds to a unique algebraic curve $R$ with
a map $ R\to \P^1$ which is ramified over the points $A,B,I$.

In general, such curves are described by subgroups of finite
index in ${\mathbb F}_2$.
Our assumption on the ramification indices at points 
$A,B$ implies that the curve $R$ corresponds to 
a subgroup of finite index in the quotient 
$\Z_2 * \Z_3$ of ${\mathbb F}_2$. The group $\Z_2 * \Z_3$
equals $\PSL(2,\Z)$. Thus local monodromy groups
over $A$-vertices can be either 1 or $\Z_3$ and over $B$ 
either 1 or $\Z_2$. 
This finishes the proof of the theorem.

\begin{coro}
The number of triangles in any $ABI$-triangulation
is equal to $2\,\deg(j_{\Gamma})$. 
Moreover,  
$2\,\deg(j_{\Gamma})=\sum_{i} v(i)$, where the summation is
over all vertices $i$ with color $I$.
\end{coro}

\begin{rem} 
A barycentric subdivision of any triangulation of
an oriented compact Riemann surface admits an 
$ABI$-coloring. We have to  color 
the initial vertices by $I$,
the vertices lying on the midpoints of  
the edges by $B$ and the vertices
inside the facets by $A$. 
This triangulation has the property that all 
vertices of color $A$ have valence 6, all vertices of color
$B$ have valence 4. However, a general 
$ABI$-triangulation, even if it does not contain 
vertices of type $A_2, B_2$, 
need not be a barycentric subdivision 
(with subsequent coloring) of
a triangulation.
\end{rem}

\begin{lem}\label{lem:jmodular}
Any arithmetic curve (an algebraic curve defined over a number field)
can be realized as a $j$-modular curve.
\end{lem}

{\em Proof.} By Belyi's theorem, for every arithmetic curve
$C$ we can find a map $f\, :\, C\to \P^1$ which is ramified over
$0,1,\infty$. Consider the triangulation $\tau(C)$ 
which is the preimage of the standard triangulation 
$\tau_{st}(\P^1)$. Consider the barycentric subdivision $\tau_b(\P^1)$ 
of $\tau(\P^1)$
induced by a map $ g \,:\,\P^1 \to \P^1$ of degree $6$ which
is ramified over $3$ points $ A,B, I$ and $g(0,1,\infty) = \infty$.
The composition $g\circ f : C \to \P^1$ is ramified only 
at $A,B,I$. The ramification indices at the
preimages of $A$ will be $3$ whereas the ramification indices 
at all the preimages of $B$ will $2$.
By Theorem 
\ref{theo:modular}, this exhibits  $C$ as a $j$-modular curve.

\begin{rem} 
The  $j$-modular structure on $C$ obtained in the
lemma corresponds to a monodromy group 
$\Gamma$ not containing
elements of finite order. Indeed, the elements of finite
order in $\Z_2* \Z_3 = \PSL(2,\Z)$ 
are conjugate either to elements of
$\Z_2$ or to elements of $\Z_3$. 
Since the ramification index over the points $A$ is 3 
everywhere the group $\Gamma$
does not contain elements conjugated to the ones in $\Z_3$.
Similar argument works for $\Z_2$.
\end{rem}

We have a closer relationship between $ABI$-triangulations
and the group $\Ga$. There is a bijection between
the set of $B_2$-vertices and 
conjugacy classes of subgroups of order 2 in $\Ga$. 
Similarly, there is a bijection between $A_2$-vertices and
conjugacy classes of subgroups of order 3 in $\Ga$. 
Finally, there is a bijection between the $I$-vertices and
conjugacy classes of  unipotent 
subgroups in $\Ga\subset \PSL(2,\Z)$. 
The generator of the unipotent subgroup is given by 
$ \left(\begin{array}{cc} 1 & v(i)/2\\ 0 & 1 
\end{array}\right)$, where $v(i)$ is the valence of the
corresponding $I$-vertex $i$.

\section{$j$-modular surfaces}
\label{sect:modular-surfaces}

In this section we study Jacobian elliptic surfaces 
such that the map $j_{\E}$ has degree 1.
Here $\tilde{\Ga}\subset \SL(2,\Z)$ is the global
monodromy group of the elliptic fibration $\E$.   
We call such surfaces $j$-modular surfaces and denote them by
$S_{\tilde{\Ga}}$.  

Consider the $j$-modular curve $M_{\Ga}$ where $\Ga$ is the 
image of $\tilde{\Ga}$ in $\PSL(2,\Z)$ under the natural projection. 
We want to solve the following problem: 
describe all surfaces $S_{\tilde{\Ga}}$ together with the
structure of a Jacobian elliptic fibration over the
$j$-modular curve $M_{\Ga}$  
such that the monodromy group $\tilde{\Ga}$ surjects onto $\Ga$.
We want to give a complete answer to this question using
the $ABI$-triangulation of $M_{\Gamma}$.

We have an exact sequence
\begin{equation}\label{eqn:ex}
0\ra \Z_2\ra \SL(2,\Z)\ra \PSL(2,\Z)\ra 1
\end{equation}
which induces a sequence 
\begin{equation}\label{eqn:exa}
0\ra \Z_2\ra \Ga'\ra \Ga\ra 1,
\end{equation}
where $\Ga'\subset\SL(2,\Z)$.

\begin{lem}
If $\Ga$ does not contain elements of order 2 then 
the exact sequence (\ref{eqn:exa}) splits. Equivalently,
the $ABI$-triangulation of $M_{\Ga}$ does not contain 
$B_2$-vertices.  
\end{lem}

{\em Proof.}
The group $\PSL(2,\Z)=\Z_2 * \Z_3$. Any subgroup 
of finite index is a finite free product of 
groups isomorphic to $\Z,\Z_2,\Z_3$. 
Assuming that $\Ga$ has no elements of order 2
we have a representation of $\Ga$ as a free product
of groups $\Z,\Z_3$. If we lift the generators of 
these free generating subgroups to elements of the
same order in $\Ga'$ we obtain a subgroup 
of  $\Ga'$ which projects isomorphically onto $\Ga$, 
in other words, a splitting of the exact sequence \ref{eqn:exa}.

\begin{rem}
All splittings differ by $\Z_2$-characters of $\Ga$ 
($H^1(\Ga, \Z_2)$) and 
the one we obtain may be not the best (this will be specified
in section \ref{sect:top-type}). Namely, the preimages of 
unipotent generators can be products of unipotent elements by 
the central element in $\SL(2,\Z)$. There may be no natural 
splitting.  
\end{rem}

We now describe a construction of an elliptic surface with 
prescribed monodromy. 
Consider the universal elliptic curve ${\cal E}^u\ra {\cal H}$
given as a quotient of $\C\times \H$ by $\Z\oplus \Z$.
The action
of $\Z\oplus \Z$ on $\C\times \lambda $ is given by 
$e_1 (z,\lambda)=(z+1,\lambda)$ and 
$e_2(z,\lambda)= (z+\lambda, \lambda)$
(here  $e_1, e_2$  are the generators of $\Z\oplus \Z$ and 
$(z,\lambda)\in \C\times \H$). The group $\SL(2,\Z)$ acts
on the universal elliptic curve, stabilizing 
the section $(0,\lambda)$.
Consider the quotient of the universal elliptic
curve  ${\cal E}^u\ra \H$ by $\Ga'$. 
We get an open surface $V'$ admitting 
a  fibration  (with a section) over the open curve 
$B' = \H/ \Ga'$, \
whose generic fiber is a smooth rational curve. The map 
${\cal E}^u\ra V'$ is ramified over a divisor $D$
which has at least two horizontal components:
$D_0$ (which is a smooth  zero-section of $V'\ra B'$) and $D_1$
which projects to $B'$ with degree 3 and is smooth 
and unramified over $B'$ in the complement of singular fibers. 
Denote by $V^o$ the open surface obtained by removing 
from $V'$ the singular fibers. The surface $V^o$ is fibered over an open curve
$B^o$ with fibers  $\P^1$. The intersection of the divisor $D$ with
each fiber consists of  exactly 4 points and $D$ is unramified over $B^o$. 

We want to define a double covering of 
$V^o$ which is ramified  
on every component of $D$. There is a correspondence between such
double coverings and special characters  $\chi\in
\Hom(\pi_1(V^o\setminus D), \Z_2)$. 
The group $\pi_1(V^o\setminus D)$ has a quotient which is 
a central $\Z_2$-extension of the free group $\pi_1(B^o)$.
This extension has a section 
(since the fibration  $V^o\ra B^o$ has a section)
and therefore it splits into a product $\Z_2\times \pi_1(B^o)$.
A character $\chi$ defining a double cover of $V^o\setminus D$  
is a character which is induced
from $\Z_2\times \pi_1(B^o)$ and which is an isomorphism on 
the central subgroup $\Z_2$ in $\Z_2\times \pi_1(B^o)$.

In other words, the restriction of $\chi$ to
the subgroup   $\pi_1(\P^1\setminus 4\,\, {\rm points})$ 
(for every fiber $\P^1$ of the fibration $V^o\ra B^o$)
is equal to the standard character of ${\mathbb F}_3$ (realized
as $\pi_1(\P^1\setminus 4 \,\,{\rm points}$) which 
sends the standard generators of ${\mathbb F}_3$ into the non-zero 
element of $\Z_2$.  

We summarize this in the diagram:
$$
\begin{array}{ccccc}
{\mathbb F}_3            & \ra   & \pi_1(V^o\setminus D) & & \\
\downarrow   &       &  \downarrow & &  \\
\Z_2           &\times & \pi_1(B^o)  & \ra &  \Z_2\\
\uparrow   &         &  \downarrow & &  \\
\Ker(\chi) & \ra     &     \Ga'    & & 
 \end{array}
$$ 
The group $\Ker(\chi)$ is a subgroup of $\Z_2 \times \pi_1(B^o)$
of index 2 and it is isomorphic to  
$\pi_1(B^o)$. This induces a map $\Ker(\chi)\ra \Ga'$.  

The character $\chi$ 
defines a double cover $W^o(\chi)$ of $V^o$.  The preimage of every 
fiber $\P^1$ of $V^o\ra B^o$ is an elliptic curve
realized as a standard double cover of this $\P^1$. 
Thus we obtain an open surface 
$W^o(\chi)$ with a structure of an 
elliptic fibration over $B^o$. All fibers are smooth.
The monodromy $\tilde{\Ga}$ of this elliptic fibration 
coincides with the image of $\Ker(\chi)$ in $\Ga'$.
If $\tilde{\Ga}$ is not equal to the whole of $\Ga'$
then the sequence \ref{eqn:exa} splits. This also means that the
character $\chi$ is induced from $\Ga'$.  

The character $\chi$ completely defines local monodromy around 
the points in $M_{\Ga}\setminus B^o$. 
Now we compactify $V^o$ keeping the structure of an elliptic fibration 
over $M_{\Ga}$ and keeping the zero section. 
Locally, in the neighborhood of $b\in M_{\Ga}$ corresponding to
singular fibers our elliptic fibration is birationally isomorphic to 
a standard fibration from the Kodaira list. The corresponding 
birational isomorphism is biregular on the complement to the
singular fiber.  The zero section is preserved under this
birational isomorphism. Now we can 
modify our initial fibration via this fiberwise transformation 
along neighborhoods of singular fibers. The resulting
surface $V$ is smooth and it admits 
a structure of a Jacobian elliptic fibration 
with the same monodromy group $\tilde{\Ga}$. 

This surface $W(\chi)$ is not unique if 
$\tilde{\Ga}$ is isomorphic to $\Ga'$. 
In this case it depends on the choice of $\chi$. 
Since we can change $\chi$ by any character of  $\pi_1(B^o)$ we 
have $2^r$ surfaces (where $r$ is the rank of $H^1(B^o,\Z)$) with
given monodromy. Removing further points from $B^o$ and 
twisting the curve by any character which is non-trivial 
at all of these points we obtain
additional moduli in our construction (of dimension equal to the
number of removed points). Thus we have moduli.  
If $\tilde{\Ga}$ is isomorphic to $\Ga$ then 
$\chi$ corresponds to the character $\Ga'\ra \Ga'/\tilde{\Ga} = \Z_2$ and 
the surface $W(\chi)$  is unique. 

\medskip 

Now we want to outline an alternative construction of
$S_{\tilde{\Ga}}$ when $\tilde{\Ga}$ does not  contain the center 
$\Z_2\subset \SL(2,\Z)$.
In this case the exact sequence \ref{eqn:exa} splits. 
Let $\tilde{\Ga}$ be any section of it.  
We realize  $\tilde{\Ga}$ as the monodromy group
of an elliptic fibration as follows:
Take the quotient $V^o\ra \H/\Ga$ of the universal elliptic curve 
${\cal  E}^u\ra \H$ by $\tilde{\Ga}$, it has the structure of a 
fibration with a section and with generic 
fibers smooth elliptic curves. 
The monodromy of this fibration (over the open curve $B=\H/\Ga$) 
is $\tilde{\Ga} \simeq \Ga$. 
Now we compactify $V^o$ keeping the 
structure of an elliptic fibration 
(over $M_{\Ga}=\overline{\H}/{\Ga}$) and the zero section as above.

\begin{rem}
It is clear that
the second construction is birationally universal
(if $\tilde{\Ga}$ does not contain the center $\Z_2$). 
Indeed, in this case, if there is a Jacobian 
elliptic fibration $V'$ with the given monodromy group
$\tilde{\Ga}$ then there is a rational fiberwise map $V\ra V'$
which is regular on the grouplike parts of $V$ and $V'$. 
\end{rem}

\section{Lifts}\label{sect:lifts}

We keep the notations of the previous sections.
Consider the diagram 
$$
\begin{array}{ccccc}
\E &  & S_{\tilde{\Ga}} & & \\
\downarrow & & \downarrow  & & \\
B  & \ra & M_{\Ga} & \ra &  \P^1 
\end{array}
$$
Let 
$B^o= B\setminus j^{-1}\{0,1,\infty\}$ and 
$M_{\Ga}^o=M_{\Ga}\setminus j_{\Ga}^{-1} \{0,1,\infty\}$  
(the points deleted from $M_{\Ga}$ 
are the $A$, $B$  and $I$ -vertices of the 
$ABI$-triangulation). 
There is a natural map 
$\pi_1(B^o)\ra \pi_1(M_{\Ga}^o)$ and a commutative diagram of
monodromy homomorphisms:
$$
\begin{array}{ccc}
\pi_1(B^o) & \longrightarrow     &      \Ga' \\
\downarrow &                     &       \downarrow \\
\pi_1(M_{\Ga}^o) & \longrightarrow &    \Ga
\end{array}
$$
and a monodromy homomorphism  $\pi_1(M_{\Ga}^o) \ra \Ga'$, 
compatible with the projection $\Ga'\ra \Ga$.

We want to compare the lifting of 
the elliptic fibration $S_{\tilde{\Ga}}$ 
to $B$ and $\E$. 
First of all we need to determine the monodromy
of a smooth relatively minimal model of the pullback  
$j^*(S_{\tilde{\Ga}})\ra B$. 
Its local monodromies
are induced by $j$ from the local monodromies of 
$S_{\tilde{\Ga}}$.
More precisely, if locally the map is given by $f(z)=z^a$
then the corresponding monodromy is $T_{z}=T^a_{f(z)}$.
Its global monodromy however can 
be different from $\tilde{\Ga}$ if $\tilde{\Ga}$ contains
the center $\Z_2$.  
Our description of the modular elliptic fibration $S_{\tilde{\Ga}}$
yields the following:

\begin{prop}\label{prop:lll}
The global monodromy  group  $\tilde{\Ga}_{B}$
of $j^*(S_{\tilde{\Ga}})\ra B$
is either isomorphic to $\tilde{\Ga}$ or
is a subgroup  of $\tilde{\Ga}$ of index 2
not containing the center $\Z_2$ (provided $\Z_2$ is a direct
summand in $\tilde{\Ga}$)
in $\tilde{\Ga}$). In the second case 
the map $B\ra M_{\Ga}$ can be decomposed into a composition
of a double covering $M_B\ra M_{\Ga}$ (which is ramified at 
some points in $ M_{\Ga}\setminus M_{\Ga}^o$) and the map 
$B\ra M_B$. The double cover $M_B\ra M_{\Ga}$ corresponds to 
the $\Z_2$-character of 
$$
\pi_1( M_{\Ga}\setminus M_{\Ga}^o)\ra \tilde{\Ga}/\Ga_B. 
$$
\end{prop}

\begin{rem}
If $g(B)=0$ then the double cover in proposition 
\ref{prop:lll} above is ramified exactly 
at two points.  
\end{rem}

\begin{prop}\label{prop:liftind}
If the local monodromies in $\E$ are
induced by $j_{\E}$ from the local monodromies of $S_{\tilde{\Ga}}$ 
and if the base $B=\P^1$ 
then $\E$ is fiberwise birationally isomorphic to 
$j^*_{\E}(S_{\tilde{\Ga}})$. 
\end{prop}

{\em Proof.} 
Consider the induced surface $j^*_{\E}(S_{\tilde{\Ga}})$.
The $j$-map coincides with the $j$-map for $\E$ and therefore
both elliptic surfaces have the same map 
$j_{\Ga}$. Since all local monodromies are the same the group 
$\tilde{\Ga}$ is mapped to 
$\Ga_B$. This is an embedding and hence an isomorphism. 
Since both the global and the local monodromies are the same we
have a fiberwise birational isomorphism between 
$\E$ and $j^*_{\E}(S_{\tilde{\Ga}})$.

In the general case, $\E$ is obtained from 
$j^*_{\E}(S_{\tilde{\Ga}})$ by 
performing an even number of twists 
corresponding to local involutions.
In particular, the two surfaces are not birational.

\begin{rem} 
If the group $\tilde{\Ga}\simeq \Ga$
then $\E$ is induced (birationally)  from $S_{\tilde{\Ga}}$. 
In particular, $h^{2,0}(\E)\ge h^{2,0}(S_{\tilde{\Ga}})$. 
\end{rem}

\begin{prop}
Assume that we have an $ABI$-triangulation $\tau$ of $\P^1$
containing vertices of type $B_2$. 
Then there exists a unique 
$\tilde{\Ga}=\tilde{\Ga}_{\E}$ corresponding 
to an $\E\ra B$ such that the corresponding $ABI$-triangulation on
$M_{\Ga}$ is isomorphic to $\tau$. 
\end{prop}

{\em Proof.} The vertices of type $B_2$ 
correspond to (conjugacy classes of) elements of $\Ga$ of order
2. The preimages of these elements in $\tilde{\G}$
are of order $4$. It follows that 
$\tilde{\Ga}$ contains a unique 
central element of order 2 and that 
$\tilde{\Ga}\in \SL(2,\Z)$ is uniquely determined by $\Ga\in \PSL(2,\Z)$. 
$\square$

If $M_{\G}$ has no vertices of type $A_2$ or $B_2$ 
then $\Ga$ is a free group. In this case, 
if $b\in B$ corresponds to a singular
fiber of $\E$ then $j_{\E}(b)$ is an $I$-vertex of the 
$ABI$-triangulations on $M_{\G}$. 
The preimage $j_{\E}^{-1}(i)$ of any $I$-vertex $i$ 
is a singular fiber of $\E$. 
Any $I$-vertex $i$ determines a (conjugacy class of a) unipotent 
element $\gamma(i)\subset\PSL(2,\Z)$  of order  $v(i)/2$. 
An element $\gamma(i)$ lifts to $\tilde{\gamma}(i)\in \SL(2,\Z)$.
The lift depends on the type of the
singular fiber at the corresponding $b(i)$; 
if the fiber $\E_{b(i)}$
is multiplicative, then $\tilde{\gamma}(i)$ is unipotent.  
Otherwise, it is $-1$ times a unipotent element.

\section{The topological type of $j$-modular surfaces}
\label{sect:top-type}

In this section we determine the  topological class of 
a $j$-modular surface $S_{\tilde{\Ga}}$ 
using the $ABI$-triangulation and the information about 
local monodromy homomorphisms.

\subsection{Degree defects}
\label{sect:degree-defects}

From now, on we assume that $B=\P^1$ (for simplicity).
Similar techniques work 
for any base $M_{\Ga}$.

Jacobian elliptic fibrations over $\P^1$ 
arise in families, defined (in Weierstrass form) as follows: 
Denote by $U_0={\mathbb A}^1 $ a chart of $\P^1$
obtained by deleting $(0:1)$ and by $U_{\infty}={\mathbb A}^1$
the chart obtained by deleting $(1:0)$.
On $U_0$ we use the coordinate $t$ and on $U_{\infty}$ the
coordinate $s=1/t$. 
Consider a hypersurface in $\P^2\times U_0$ given by
$$
zy^2=x^3+p_0(t)xz^2+q_0(t)z^3
$$
where $p$ (resp. $q$) is a polynomial of degree $4r$ (resp. $6r$). 
In $U_{\infty}$ the equation is
similar, with $p_{\infty}(s)=p_0(1/s)s^{4r}$ and $q_{\infty}(s)=
q_0(1/s)s^{6r}$. We get elliptic fibrations over $U_0, U_{\infty}$
which we can glue to an elliptic surface ${\cal E}\ra B$.
The $j$-function (on $U_0$) is given by
$p_0(t)^3/(4p_0(t)^3+27q_0(t)^2)$.  
The obtained fibration can be singular in 
fibers corresponding to 
$b\in B$ where $4p_0(t)^3+27q_0(t)^2=0$ and the singularities
can be resolved by a sequence of blow-ups. 
The outcome is a (unique) smooth relatively minimal 
Jacobian elliptic fibration. 
Thus we get a family 
${\cal F}_r$ of such elliptic fibrations.  
Notice  that $12 r=\chi({\cal O}_{\cal E})$. 
Conversely, a simply connected, compact, minimal 
Jacobian elliptic fibration
with $\chi({\cal O}_{\cal E})=12 r$ belongs to ${\cal F}_r$. 
The family ${\cal F}_r$ is 
parametrized by
the coefficients of $p_0,q_0$ (subject to certain constrains)
- it is a smooth irreducible variety.
Every  Jacobian elliptic fibration is birational to 
a minimal elliptic fibration and the $j$-map for both
fibrations is the same.

\

The generic degree of the $j$-map in the family 
${\cal F}_r$ is $12r$. 
However, the presence of fibers of non-multiplicative
type diminishes the degree of $j$. 
We define the degree defect as
$$
\df(\E):=12 r - \deg(j).
$$ 
This degree defect results from 
possible common roots of $p^3(t)$ and 
$4p(t)^3+27q^2(t)$ in the formula $j(t)$ for $t$ 
corresponding to singular fibers 
and therefore, $\df(\E)$ is a sum of local contributions from 
singular fibers
of $\E$. We denote by $\df(\E_b)$  (for $b\in B^s$) 
these local contributions. 

\begin{prop} \label{prop:defects}
Let $\E\ra B$ be an elliptic fibration over $\P^1$
and $\E_b$ be a singular fiber of $\E$. 
\begin{enumerate}
\item If $\E_b$ is of type {\rm II}, {\rm III} or {\rm IV }
then the local contribution $\df(\E_b)$ 
is at least $2,3$ or $4$, respectively.
\item  If $\E_b$ is a quotient of 
a fiber of type {\rm II}, {\rm III}, {\rm IV} 
by the action of the birational involution ($x\ra -x$)
(these fibers are denoted by ${\rm II}^*, {\rm III}^*, {\rm IV}^*$)
then $\df({\cal E}_b)$ 
is at least 8,9 or 10, respectively. 
\item If  ${\cal E}_b$ is a singular fiber which is 
a quotient of a smooth or multiplicative fiber by the involution 
(these fibers are denoted by ${\rm I}_0^*, 
{\rm I}_n^*$), then $\df(\E_b)$ is at least 6. 
\end{enumerate}
\end{prop}

{\em Proof.} Local computation, see \cite{shafarevich}, p. 171.

For $\E=S_{\tilde{\Ga}}\ra M_{\Ga}$ 
we can translate the topological information into
the combinatorics of the $ABI$-triangulation of $M_{\Ga}$. 
Notice that for $S_{\tilde{\Ga}}$ the degree estimates of 
Proposition \ref{prop:defects} are sharp. 
We have the following

\begin{prop}
Assume that $M_{\Ga}$ does not contain $A_2$ or $B_2$-vertices
and that the generators of all monodromy groups $T_b$ are unipotent.
Then any $j$-modular surface $S_{\tilde{\Ga}}$ over
$M_{\Ga}$ belongs to ${\cal F}_r$, with $24r-12n$ 
equal to the number $|\tau_2(\Ga)|$ of 
open two cells of the $ABI$-triangulation 
of $M_{\Ga}$, where $n$ is the number of fibers twisted by
the involution (see Section \ref{sect:modular-surfaces} for the
construction of $S_{\tilde{\Ga}}$).  
\end{prop}

{\em Proof.} 
Indeed, we can compute the Euler characteristic of the semistable
fibration $S_{\tilde{\Ga}}$. In this case,
the  singular fibers of $S_{\tilde{\Ga}}\ra M_{\Ga}$ 
lie exactly over
the $I$-vertices of the $ABI$-triangulation of 
$M_{\Ga}$ and these fibers
are all of multiplicative type 
${\rm I}_n$. Here $n$ is equal to $v(i)/2$ of 
the corresponding $I$-vertex $i$. 
The Euler characteristic is given as $\sum n_i$ over the set
of $I$-vertices $i$. The set of all $ABI$-triangles
is a disjoint union of the sets of 
triangles having one $I$-vertex in common. 
The number of triangles in the latter 
set is equal to the valence of the 
corresponding $I$-vertex. 
Since the contribution to the Euler characteristic 
of the singular fiber over this $I$-vertex 
$i$ is $v(i)/2$ we get the result. 

\

This proposition not only shows how to compute the class of 
a $j$-modular surface in the ideal situation when all singular
fibers are of multiplicative type but also demonstrates non-trivial
combinatorial restrictions on the $ABI$-triangulations corresponding 
to such $j$-modular fibrations. Notice that in this case
the degree of the map $j_{\Ga}$ is equal to $12r-6n$.

\

The presence of $A_2, B_2$-vertices on $M_{\Ga}$ 
diminishes the degree of $j_{\Ga}$ (which equals
half the number of triangles in the $ABI$-triangulation). 
It will be convenient for us to use the combinatorial analogs
of degree defects to estimate from below the change of the
degree. 
For any $ABI$-triangulation of $M_{\Ga}$
we define the combinatorial degree defect
as follows: 
$$
\cdf(\Ga)= 2 a_2 + 3 b_2,
$$  
where $a_2$ (resp. $b_2$) is the number of 
$A_2$ (resp. $B_2$) vertices
in the $ABI$-triangulation on $M_{\Ga}$. 
Denote by $\ET(M_{\Ga})$ the number of ``effective triangles''
of the $ABI$-triangulation corresponding $\Ga$. 
By definition, 
$$
\ET(M_{\Ga}) = |\tau_2(\Ga)|+ 2 \cdf(\Ga),
$$
(where $\tau_2(\Ga)$ is the number of open 2-cells
in the $ABI$-triangulation).
Notice that Proposition \ref{prop:defects} and

\begin{lem}
We have
$$
\df(S_{\tilde{\Ga}})= \cdf(\Ga)+6n,
$$
where $n$ is the number of fibers of $S_{\tilde{\Ga}}$ 
twisted by the involution. 
\end{lem}

{\em Proof.} 
The result follows from Proposition 
\ref{prop:defects} and the fact that
$A_2$, $B_2$-vertices in $M_{\Ga}$ correspond 
to singular fibers of $S_{\tilde{\Ga}}\ra M_{\Ga}$ of the 
types listed in that proposition.

\begin{coro}
If $S_{\tilde{\Ga}}\in {\cal F}_r$ then
$12r\ge \ET(\Ga)$. 
\end{coro}

The actual degree defect of $S_{\tilde{\Ga}}$ depends
not only on the $ABI$-triangulation  (which determines 
$\Ga$) but also on the choice of a lifting of 
local monodromies to $\tilde{\Ga}$. This 
leads us to: 

\begin{dfn}
We shall say that $S_{\tilde{\Ga}}$ is {\em minimal} if
for every singular fiber the contribution to 
the local degree defect corresponding to 
lifting (to $\tilde{\Ga}$) of local monodromy (from $\Ga$)
is minimal. 
\end{dfn}

Now we give combinatorial criterium for the existence
of minimal $S_{\tilde{\Ga}}$.

\begin{theo}\label{theo:twist} 
If $M_{\Gamma} = \P^1 $ then 
a minimal lifting exists iff 
$\ET(M_{\Gamma})$ is divisible by $24$. Since $\ET(M_{\Gamma})$
is always divisible by $12$ there is always a lifting
with at most one non-minimal local monodromy element.
\end{theo}

{\em Proof.}
Every vertex $v$ of the $ABI$-triangulation 
determines a standard element $T_v'$ in $\PSL(2,\Z)$. 
In our construction of $S_{\tilde{\Ga}}$
we had a choice of two possible local monodromies
for each fiber - these correspond to two
possible lifts  $T_v$ of $T_v'$ to $\SL(2,\Z)$. 
The difference in corresponding degree defects is 6. 
One of the liftings is minimal, with respect to $\df({\cal E}_b)$
(for the fiber ${\cal E}_b$). 

It remains to find out when it is possible to 
choose these minimal liftings
compatibly for all fibers. 
Compatibility is equivalent to 
$\prod T_v=1$ in $\SL(2,\Z)$.  
Since we know that $\prod T_v' =1 $ in $\PSL(2,\Z)$ the only
possibility is that $\prod T_v $ is 1 or the central element
$c\in \SL(2,\Z)$.
To determine when the product is equal to $c$ we 
use the existence of a standard lifting of
Dehn twists into the group $\tilde{\SL}(2,\Z)$. Here
we denote by $\tilde{SL}(2,\Z)$
the preimage of 
$\SL(2,\Z)$ in the universal cover of $\SL(2,\R)$.
More precisely, local 
monodromies $T_v$ can be represented as finite products
of right Dehn twists.  

\begin{lem}
Denote by $d_v$ the number of Dehn twists representing $T_v$.
The sum $\sum d_v$ is always divisible by 6. 
If $\sum d_v$ is divisible by 12 then the product
$\prod T_v =1$ in $\SL(2,\Z)$.  
\end{lem}

{\em Proof.}
A standard right Dehn twist is 
conjugated to the element 
$\left(\begin{array}{cc} 1 & 1 \\ 0 & 1 \end{array}\right) 
$ in $\SL(2,\Z)$. Each such element has a standard lifting 
into $\tilde{\SL}(2,\Z)$. 
The group $\tilde{\SL}(2,\Z)$ is a braid group generated
by standard right Dehn twists $a,b$  with one
braid relation $aba=bab$. 
Therefore, we have a  (degree) homomorphism 
$\chi \, :\, \tilde{\SL}(2,\Z)\ra \Z$. 
The image
of all Dehn twists in $\Z$ is equal to 1. 
The generator $\tilde{c}=(aba)^2$ 
of the center of $\tilde{\SL}(2,\Z)$ projects
into the center $\Z_2$ of $\SL(2,\Z)$.
We have $\chi(\tilde{c})=6$. 

Denote by $\tilde{T}_v$ 
the liftings of local monodromies into 
$\tilde{\SL}(2,\Z)$. 
Thus we have a well defined $\chi(\tilde{T}_v)\in \Z$.
Assuming that the product of local monodromies
$\prod T_v =1 $ in $\SL(2,\Z)$ we see that
$\chi(\prod \tilde{T}_v)$ has to be divisible by 12.   
Therefore, the sum of the number of Dehn twists representing
local monodromies has to be divisible by 12. 

Since $\chi(\tilde{c})$ is equal to 6, $\sum d_v$ 
is always divisible by 6.

To conclude the proof of Theorem \ref{theo:twist}
it suffices to observe that 
the number of Dehn twists in the decomposition 
of minimal local monodromies of finite order is 
equal to the degree defect of the corresponding singular
fiber and that the number of 
Dehn twists representing the unipotent monodromies
is equal to $ 1/2 v(i)$  of the 
corresponding $I$-vertex of the $ABI$-triangulation.

Therefore, the lemma implies that 
if $\ET(\Ga)$ is divisible by 24 then the product
of local monodromies is equal to 1. 
If it is divisible by 12 we can twist some fiber
to obtain the relation $\prod T_v=1$, increasing the
degree defect by 6.

\section{Combinatorics}

\subsection{Divisibility by 12}
We start with an $ABI$-triangulation on $M_{\Ga}=\P^1$ and remove
the $I$-vertices, together with all the connections to the
$A$ and $B$-vertices. We obtain an 
$AB$-graph which we draw on the plane.  It might  look as follows:

\vskip 0,5cm
\hskip 2cm\includegraphics{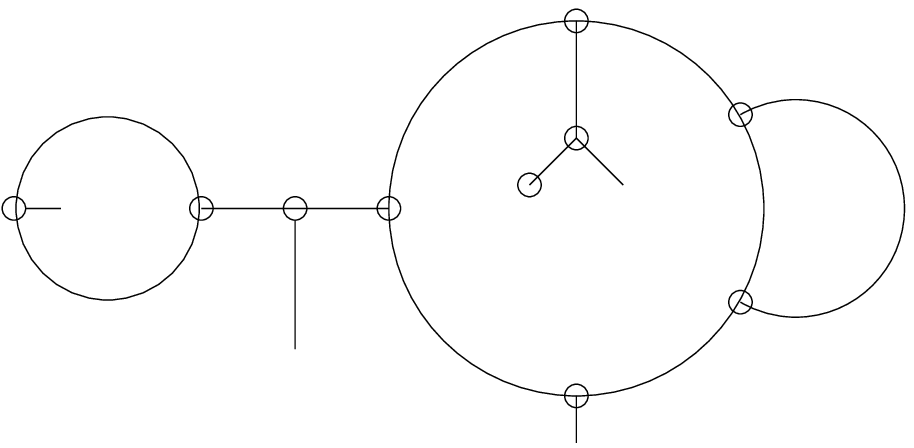}
\vskip 0,5cm

Here we use a small circle to indicate an $A$-vertex.  
The $B$-vertices
are placed on the edges between two $A$-vertices. A ``loose'' end
represents a $B$-vertex. 

Every $AB$-graph can be simpliefied as follows:
clipp off all trees together
with the vertex where they originate. 
The outcome is a connected graph 
without ends and with only $A_6$-vertices. 
We can 
think of the remaining graph as coming from a ``generalized''
triangulation of $\P^1$. These are elementary objects
with $\ET=6a_6$ 
(where $a_6$ is the number of $a_6$-vertices).
$\ET$ of the initial graph is sum of $\ET$ from the 
trees + $\ET$ of the elementary graph obtained.

\begin{lem} 
Consider $M_{\Ga}=\P^1$ with its $ABI$-triangulation. Then 
$\ET(\Ga)$ is divisible by 12.
\end{lem}

{\em Proof.} 
Recall that $\ET= 2|AB|-{\rm edges} + $ contributions from 
vertices. 
First, we remove all $B_2$-vertices and the adjacent
edges. The value for $\ET(\Ga)$ 
changes by a multiple of 12 (simple check). 
Next we pick an $A_6$-vertex and disconnect the $AB$-graph, 
removing two of the edges adjacent to it. 
We obtain an $A_2$-vertex. 
The new $\ET_{\rm new}= \ET_{\rm old}-8+8+12$. 

\vskip 0,5cm
\hskip 2cm\includegraphics{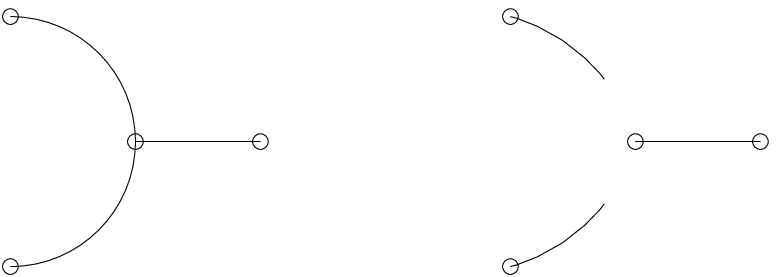}
\vskip 0,5cm

Now again we remove all $B_4$-vertices. This way we continue
until there are no $A_6$-vertices. 
The outcome is a collection of simple chains of the type:

\vskip 0,5cm
\hskip 4,5cm\includegraphics{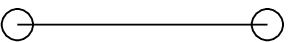}
\vskip 0,5cm

The contribution from such chains is 
$12$.

\subsection{Graphs without loops}

It is easy to compute $\ET(\Ga)$ if the corresponding $AB$-graph
has no loops. Indeed, such graphs 
are represented as follows: 
take any tree with only triple ramifications and 
mark arbitrarily some ends by $A_2$.

\begin{lem}
If the $AB$-graph has $k$ ends  then the number
of $A_6$-vertices is $k-2$ and 
$$
\ET(\Ga)=6k+6(k-2).
$$
\end{lem}

{\em Proof.} 
Every ramification in the graph has to be an $A_6$-vertex.
The ends can be either $A_2$ or $B_2$-vertices - 
the corresponding contributions to $\ET$ are either 4 or 6.

\begin{rem}
There are very few  graphs without loops and with $\ET(\Ga)\le 48$. 
The number of ends of the tree is $\le 5$.
If the number of ends is 2, 3 or 4 then there is only 
one tree, if it is 5 then there are only two trees.
\end{rem}

\subsection{Graphs with loops and small $\ET$}

First we list graphs with $\ET(\Ga)=12$:

\vskip 0,5cm
\hskip 0,5cm\includegraphics{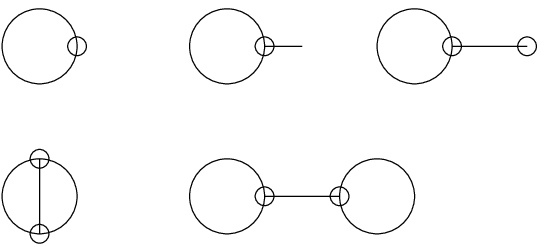}
\vskip 0,5cm

We will call an $AB$-graph saturated if all the $A$-vertices
are $A_6$-vertices. Saturated graphs can be considered as
arising from generalized triangulations of $\P^1$. 
An arbitrary graph can be obtained from a saturated graph
by addition of trees. 
It is easy to control the change of $\ET$ under this
basic operation. 
We add the tree as follows: pick a new point on one of the edges and
make it to an $A_6$-vertex with the tree attached. The $\ET$ is
the sum of $\ET$ of initial graph plus $\ET$ of the tree. 

\

To conclude, we list saturated graphs with $\ET(\Ga )=24$:
\vskip 0,5cm
\hskip 0,5cm\includegraphics{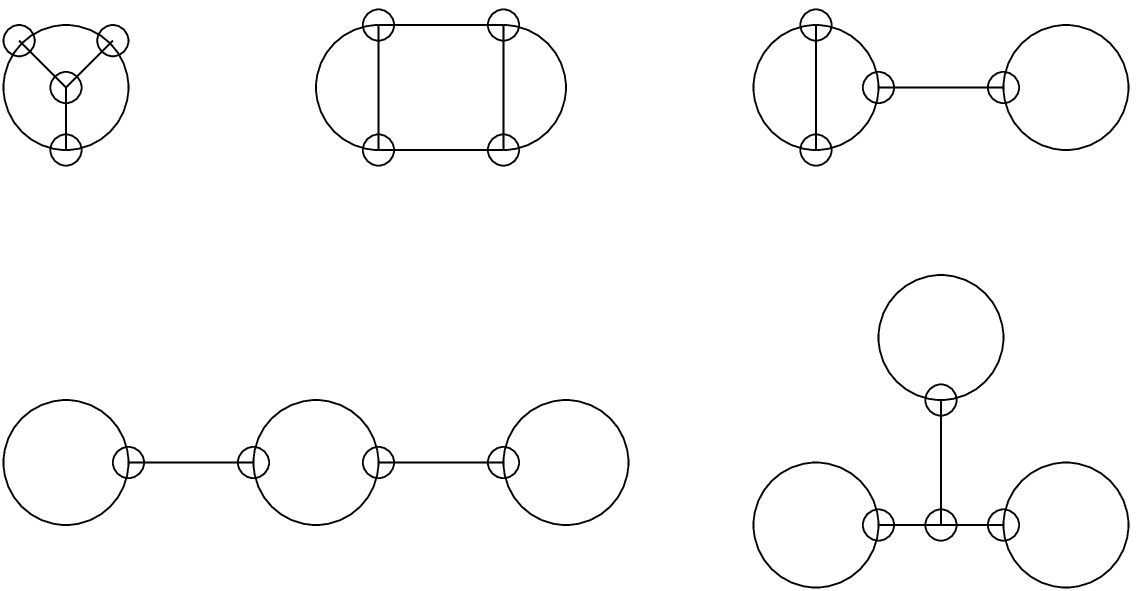}
\vskip 0,5cm

\vskip 0,5cm
\vskip 0,5cm

\vskip 0,5cm
\vskip 0,5cm

\end{document}